\documentstyle[12pt]{article}
\newcommand{\Hom}{\mbox{\rm Hom}}

\newcommand{\Zp}{\mbox{${\bf Z}_{p}$}}
\newcommand{\Z}{\mbox{${\bf Z}$}}
\newcommand{\Q}{\mbox{${\bf Q}$}}

\newcommand{\Fp}{{\bf F}_{p}}

\newtheorem{lemma}{Lemma}
\newtheorem{formula}{Formula}

\newenvironment{proof}{
\medskip
\noindent {\em Proof\/}. }{\nolinebreak[2]\rule{.35em}{.9em}\medskip}

\begin{document}
\noindent
{\sl This is the body of a letter that I (Ken Ribet)
wrote to J.-F. Mestre in early
November, 1987.  I did not publish the theorems in this
letter because I expected that they would be incorporated
into an article by Mestre and J.~Oesterl\'e on the ``graph method''
for calculating eigenvalues of Hecke operators on spaces of
cusp forms.  My understanding is that the formulas in this letter
have not appeared
in print and that they retain some interest.
I am happy to make them available at this time (May, 2001).
My formulas have been generalized to some extent by William A. Stein,
who calculates component groups of optimal abelian variety quotients
of Jacobians of modular curves.}
\bigskip

I return, first, to the notations you had on the board: $p$ is a prime
number, $M$ is prime to $p$, $N$ is the product $pM$.  We suppose given
an elliptic curve $E$, plus a map $\pi :X_{o}(N) \rightarrow E$,
over \Q .  We assume that the induced map on Albanese varieties
$\pi _{*}:J_{o}(N) \rightarrow E$ has connected kernel.  Finally, we
assume that $E$ has bad (hence multiplicative) reduction at $p$.
This situation occurs, in particular, if $E$ has conductor $N$
and is given as a ``strong (Taniyama-)Weil curve.''

Next, we have some notation relative to the N\'eron models of $E$
and $J = J_{o}(N)$ at $p$.  I will let $X$ be the character group of
the torus attached to $J_{\Fp}$, so that $X$ is the group of degree-0
divisors on the set $\cal C$ of isomorphism classes of supersingular
elliptic curves over ${\overline{{\bf F}}}_{p}$, which are enhanced
by $\Gamma _{o}(M)$ structures.  A typical
element of $X$ is a sum $\sum \lambda _{C} (C)$, where the $\lambda _{C}$ are
integers which sum to 0.  
Let $X_{E}$ be the character group associated with $E$.
We chose a generator of this cyclic group.  The natural map
$\pi ^{*}:X_{E} \rightarrow X$ then sends this generator to a certain element
$\sum \lambda _{C} (C)$ of $X$.  We can identify $X_{E}$ with the group
of multiples of this sum.  There is a second natural map 
$\pi _{*}: X \rightarrow X_{E}$.  The endomorphism of $X_{E}$
gotten by following $\pi ^{*}$ by $\pi _{*}$ is multiplication by
$n = {\rm deg}(\pi )$.

There are natural symmetric
 \Z -valued pairings $\langle\: ,\:\rangle$ and
$\langle\: ,\:\rangle _{E}$ on $X$ and $X_{E}$.  These are the
{\em monodromy pairings} considered by Grothendieck in SGA7.  The component
groups $\Phi$ and $\Psi$ attached to $J$ and $E$ are then the cokernels
of the maps \[X \rightarrow \Hom (X,\Z )\] and 
\[X_{E} \rightarrow \Hom (X_{E},\Z )\] deduced from these pairings.
The map \[\pi _{*}:\Phi \rightarrow \Psi\] that $\pi$ induces on
component groups can be seen in the commutative diagram
\begin{equation}\begin{array}{cclccclcl}
0&\rightarrow&X&\rightarrow&\Hom (X,\Z )&\rightarrow&\Phi&\rightarrow&0\\
\null &\null &\downarrow  \pi _{*} &\null &\downarrow  \beta&
\null&\downarrow  \pi _{*}&\null&\null \\ 
0&\rightarrow&X_{E}&
\rightarrow&\Hom (X_{E},\Z )&\rightarrow&\Psi&\rightarrow&0,\\
\end{array}\end{equation}
in which $\beta$ is Hom$(\pi ^{*},\Z )$.

\begin{lemma}
The map $\pi _{*}:X \rightarrow X_{E}$ is surjective.
\end{lemma}
\begin{proof}
The map (in characteristic 0) $\pi ^{*}:E \rightarrow J$ is
injective, because its transpose $\pi _{*}$ has connected kernel,
by assumption.  In particular, the map $E[\ell ] \rightarrow J[\ell ]$
is injective, for all prime numbers $\ell$.  The map which $\pi ^{*}$
induces on tori \[E^{\rm toric} \rightarrow J^{\rm toric}\] is then
injective (in the appropriate sense), as one sees by lifting these
tori to \Zp .  This injectivity corresponds to the
surjectivity asserted by the lemma.
\end{proof}

\begin{formula}
The cokernel of $\pi _{*}:\Phi \rightarrow \Psi$ has order
$\gcd (\ldots ,\lambda _{C},\ldots )$.
\end{formula}

\begin{proof}
The cokernel in question is
\[
\Hom (X_{E},\Z ) /\{ X_{E} + \beta (\Hom (X,\Z            )  \}.
\]
By the lemma, we have $X_{E} \subseteq \beta (\Hom (X,\Z )$, since
$\Hom (X,\Z )$ contains $X$.  Therefore the cokernel is simply
\[
\Hom (X_{E},\Z ) /\beta (\Hom (X,\Z ).
\]
After some fooling with  Ext groups, we see that this quotient
map be identified with the torsion subgroup of $X/X_{E}$
\footnote{Actually, this claim is wrong; I think that you
get instead the {\em dual\/} of the indicated torsion group}.
However, the gcd in the formula is 
the order of the torsion subgroup of $X/X_{E}$.
\end{proof}

\begin{formula}
We have
\begin{eqnarray}
{\rm card}(\Psi )  = & \gcd &(\lambda _{C} w_{C} - \lambda _{C'} w_{C'}),
\nonumber \\
\null & \null ^{C,C'} &\null \nonumber
\end{eqnarray}
where $w_{C}$ is the integer {\rm Aut}$(C)/2$, as in your notation.
\end{formula}
\begin{proof}
Consider the map
\begin{equation}
 \gamma :X \rightarrow \Hom (X_{E}, \Z ) 
\end{equation}
gotten by composing \[\pi _{*}: X \rightarrow X_{E}\] and
the map \[X_{E} \rightarrow \Hom (X_{E}, \Z ) \] which is 
$\langle\: ,\:\rangle _{E}$.  By the lemma, its cokernel is $\Psi$.
To obtain the formula, we notice that the image of 
$\gamma$ is generated by the elements $\gamma (D-D')$
for $D,D' \in {\cal C}$.  For $D$ and $D'$ in $\cal C$, $\gamma (D-D')$
maps the generator $\sum \lambda _{C} (C)$ of $X_{E}$ to
$\lambda _{D} w_{D} - \lambda _{D'} w_{D'}$, since $\langle\: ,\:\rangle$
is deduced by restriction from the pairing
 $C,C' \mapsto \delta _{C,C'}w_{C}$ on $\Z ^{{\cal C}}$.
It should be remarked here that we are using the compatibility between
$\langle\: ,\:\rangle _{E}$ and $\langle\: ,\:\rangle $ which is summarized
by the relation
\begin{equation}
\langle \xi , \pi _{*}(\eta )\rangle _{E}  =
  \langle \pi ^{*}(\xi ) ,\eta \rangle 
\end{equation}
for $\xi \in X_{E}$ and $\eta \in X$.
\end{proof}

\begin{formula}
We have
\[
n\: {\rm card}(\Psi ) =  \sum \lambda _{C}^{2} w_{C} .
\]
\end{formula}
\begin{proof}
Let $g$ be our chosen generator of $X_{E}$.  Then 
\[{\rm card}(\Psi ) = \langle g , g \rangle _{E},
\] so that 

 \[ n\:{\rm card}(\Psi ) = \langle \pi _{*}\pi ^{*}g , g \rangle _{E}
= \langle \pi ^{*}g , \pi ^{*}g \rangle . \]
The quantity $\langle \pi ^{*}g , \pi ^{*}g \rangle$ visibly coincides
with the sum $\sum \lambda _{C}^{2} w_{C}$.
\end{proof}
\par
\bigskip
{\sf Ken RIBET}
\end{document}